\documentclass[12pt]{amsart}
\usepackage{amsxtra}
\usepackage{amssymb}
\newcommand{\cleqn}{\setcounter{equation}{0}}
\newcommand{\clth}{\setcounter{theorem}{0}}
\newcommand {\sectionnew}[1]{\section{#1}\cleqn\clth}
\newcommand {\subsec} {\subsection}
\theoremstyle{plain}
\newtheorem{theorem}{Theorem}[section]

\newtheorem{definition-theorem}[theorem]{Definition-Theorem}

\theoremstyle{definition}
\newtheorem{definition}[theorem]{Definition}

\newtheorem{example}[theorem]{Example}
\newtheorem{examples}[theorem]{Example}

\newtheorem*{remark*}{Remark}

\newtheorem{observation}[theorem]{Observation}

\newtheorem{problems}[theorem]{Open problems}

\numberwithin{equation}{theorem}

\input xy \xyoption{matrix} \xyoption{arrow}
          \xyoption{curve}  \xyoption{frame}

\def\edge{\ar@{-}}
\def\dttdar{\ar@{.>}}
\def\pbl{\save+<0ex,0ex> \drop{\bullet} \restore}
\def\drbl{\save+<0ex,-2ex> \drop{\bullet} \restore}
\def\horizpool#1{  \save[0,0]+(0,3);[0,0]+(0,3) **\crv{~*=<2.5pt>{.} [0,#1]+(0,3) &[0,#1]+(3,3) &[0,#1]+(3,0) &[0,#1]+(3,-3) &[0,#1]+(0,-3) &[0,0]+(0,-3) &[0,0]+(-3,-3) &[0,0]+(-3,0) &[0,0]+(-3,3)}\restore }
\def\dropdown#1{\save+<0ex,-3ex> \drop{#1} \restore}

\def\AA{\mathbb{A}}
\newcommand\CC{\mathbb{C}}

\def\SS{\mathbb{S}}
\newcommand\ZZ{\mathbb{Z}}
\newcommand\NN{\mathbb{N}}
\newcommand\QQ{\mathbb{Q}}

\newcommand\A{\mathcal{A}} 
 
\newcommand\C{\mathcal{C}}
\newcommand\D{\mathcal{D}}
\newcommand\E{\mathcal{E}}
\newcommand\F{\mathcal{F}}

\def\L{\mathcal{L}}

\def\P{\mathcal{P}}

\newcommand\fancyS{\mathcal{S}}

\newcommand\la{{\Lambda}}
\newcommand\latilde{\widetilde{\la}}

\def\lamod{\Lambda\mbox{\rm-mod}}
\def\lalmod{\Lambda_L\mbox{\rm-mod}}
\def\deltamod{\Delta\mbox{\rm-mod}}

\def\adm{\text{\bf adm}(\lamod)}

\def\etil{\tilde{e}}

\def\Qtilde{\widetilde{Q}}
\def\bd{{\mathbf d}}

\DeclareMathOperator{\Filt}{\mathbf{Filt}}
\DeclareMathOperator{\GL}{GL} 
\DeclareMathOperator{\End}{End}

\DeclareMathOperator \Hom{{Hom}}

\DeclareMathOperator\soc{{soc}}

\DeclareMathOperator\gldim{{gl.dim}}

\DeclareMathOperator\Ker{{Ker}}
\def\Im{\mbox{\rm Im}}
\DeclareMathOperator\rep{{\mathbf{Rep}}}

\DeclareMathOperator\add{{add}}
\DeclareMathOperator\Ext{{Ext}}

\DeclareMathOperator\pdim{{p\, dim}}

\DeclareMathOperator\udim{{\underline{\dim}}}
\DeclareMathOperator\op{{op}}
\DeclareMathOperator\ttop{{top}}

\def\pinf{{\P^{< \infty}}}

\def\eps{\varepsilon}
\newcommand\onetil{\widetilde{1}}
\newcommand\twotil{\widetilde{2}}
\newcommand\threetil{\widetilde{3}}
\newcommand\fourtil{\widetilde{4}}

\begin{document}

\title{Finite Dimensional Representations of Quivers with Oriented Cycles}

\author[Goodearl]{\ K.R.  Goodearl}
\address{Department of Mathematics, University of California, Santa Barbara, CA 93106, USA}
\email{goodearl@math.ucsb.edu}

\author[Huisgen-Zimmermann]{\ B. Huisgen-Zimmermann}
\address{Department of Mathematics, University of California, Santa Barbara, CA 93106, USA}
\email{birge@math.ucsb.edu}

\begin{abstract}  Let $K$ be a field, $Q$ a quiver, and $\A$ the ideal of the path algebra $KQ$ that is generated by the arrows of $Q$. We present old and new results about the representation theories of the truncations $KQ/ \A^L$, $L \in \NN$, tracking their development as $L$ goes to infinity. The goal is to gain a better understanding of the category of those finite dimensional $KQ$-modules which arise as finitely generated modules over admissible quotients of $KQ$.   
\end{abstract}

\dedicatory{Dedicated to Laszlo Fuchs on the occasion of his 100th birthday}



\maketitle

\sectionnew{\bf Introduction and notation}

\noindent This is an overview of recent work on the finite dimensional representation theory of infinite dimensional path algebras $\la = KQ$, i.e., of path algebras of quivers $Q$ with oriented cycles.  Throughout, the base field $K$ is algebraically closed. Some of the results we present here concerning strong tilts have not yet been published, so this article also serves as a first announcement.  We do not include proofs here, but illustrate the theory with numerous examples.

Due to Gabriel \cite{Gab}, all basic finite dimensional algebras are quotients of path algebras, whence path algebras hold a foundational position within finite dimensional representation theory.  As a consequence, they have been thoroughly explored over the past forty-five years. They are known to have multiple assets, one of them being $\gldim \la  \le 1$, another being the role played by the Euler form of $Q$; indeed this form connects the dimensions of Ext-spaces of finite dimensional $\la$-modules to those of their Hom-spaces. We refer to \cite{AsSiSk, ARS, CB.reps.quivers, CB.morereps.quivers} for background.  Several of these assets are restricted to path algebras of acyclic quivers, one of the exceptions being the geometry of the varieties $\rep_\bd(\la)$ parametrizing the modules of fixed 
dimension vector $\bd$, and the generic decompositions of the modules they encode. 

Asking for a full understanding of the category $\lamod$ of finite dimensional left $\la$-modules is not an attainable goal; after all, the free algebras in any finite number of variables are among the path algebras.  In a first step, we narrow the problem to targeting only those finite dimensional $\la$-modules which are annihilated by some power of the ideal $\A := \langle \text{all arrows of}\  Q \rangle \subseteq KQ$.  We call these modules \emph{admissible}, and focus on the full subcategory $\adm$ of $\lamod$ consisting of the admissible modules.  In light of the fact that every basic finite dimensional $K$-algebra is isomorphic to a quotient $KQ/I$, where $I$ is an admissible ideal of $KQ$ and as such contains a power of $\A$ \cite[\S4.3]{Gab}, the object class of $\adm$ is the union of the classes of finitely generated left $\Delta$-modules, where $\Delta$ traces all basic finite dimensional $K$-algebras with Gabriel quiver $Q$. 

We tackle the category of admissible $\la$-modules by first studying the finitely generated modules over the truncations $$\la_L := KQ/ \langle \text{all paths of length}\ L \rangle$$
for arbitrary positive integers $L$. That is, we explore $\adm$ along the chain of full subcategories 
$$\la_1\mbox{\rm-mod} \subseteq \la_2\mbox{\rm-mod} \subseteq \la_3\mbox{\rm-mod} \subseteq  \cdots \subseteq \adm,$$ 
the union of the object classes of the $\lalmod$ being the class of admissible modules.
In a next step, we then trace the development of the resulting pictures as $L \rightarrow \infty$.  Evidently, each truncation depends only on $Q$ and $L$, and the category $\adm$ can be seen as a ``direct limit" of the directed family of categories $(\la_L\text{-mod})_{L \in \NN}$. As  becomes clear in Sections 2 and 3, the individual $\la_L\text{-mod}$ retain several of the assets of $\lamod$, which allow for effective tracking of their representation theories, relative to each other, as $L$ grows.  

More specifically: Via the family $\bigl(\la_L\bigr)_{L \in \NN}$, we investigate the category $\adm$ from the following angles:
\smallskip

{\bf (A)} {\bf by way of the parametrizing varieties for its objects}; 
\smallskip

{\bf (B)} {\bf through the lens of tilted algebras}; 
\smallskip

{\bf (C)} {\bf in terms of representation type}.
\smallskip

\noindent  Results regarding the first two points are at an advanced level, whereas those aimed at the third are still rudimentary.  We include a short section on this topic at the end, in hopes of stimulating further work on the subject.  
\smallskip

Under {\bf (A)} we describe the irreducible components of the classical varieties $\rep_{\bd}(\la_L)$ for any dimension vector $\bd$ and provide information on the generic behavior of the modules encoded by these components; in doing so, we modestly supplement previously available facts (cf.~\cite{closures}). Recall that the modules $M_x$ corresponding to the points $x$ of an irreducible component $\C$ of $\rep_{\bd}(\la_L)$ are said to \emph{generically} have a certain property if there exists a dense open subset $U$ of $\C$ such that all modules $M_x$ with $x \in U$ share this property. 

Our classification of the irreducible components of $\rep_{\bd}(\la_L)$ is purely in terms of module invariants and allows for an effective listing from $Q$ and $L$.  Moreover, for $L < L'$, all containments $\C_i \subset \D_j$, where the $\C_i$ trace the irreducible components of $\rep_\bd(\Lambda_L)$ and the $\D_j$ those of $\rep_\bd(\Lambda_{L'})$, can be determined.
\smallskip

Under {\bf (B)} we exploit the known fact that, for any truncated path algebra $\la_L$, the category $\pinf(\la_L\text{-mod})$ of finitely generated $\la_L$-modules with finite projective dimension is contravariantly finite in $\lalmod$ \cite{DHZ}, whence by \cite{AuRe}, $\lalmod$ contains a strong tilting module (see Section 3 for definitions and background).  The corresponding tilted algebra $\latilde_L$ is particularly tightly tied to $\la_L$; it captures the finite dimensional representation theory of $\la_L$, and consequently that of $\la$, from a shifted angle.  
\smallskip

Under {\bf (C)} we address the problem of recognizing the representation type of the algebras $\la_L$ based on $Q$ and $L$, as well as the correlations between the representation types of the truncations $\la_L$ and those of their strong tilts $\latilde_L$. The problems we propose are wide open for $L \ge 3$; however, some preliminary results are available for $L = 2$.

\subsec*{Further notation and terminology}
Throughout, $\la = KQ$, where $K$ is an algebraically closed field. We denote by $e_1, \dots, e_n$  the vertices of $Q$ and systematically identify them with a full set of orthogonal primitive idempotents of $\la$.  Our convention for the multiplication of paths in $\la$ is to concatenate from right to left, i.e., $pq$ stands for ``$p$ after $q$" if $\text{start}(p) = \text{end}(q)$, while $pq=0$ otherwise.  We use the same identifications and conventions for the truncations $\la_L$; the relevant identifications are always clear from the context.  Moreover, for $i \le n$, we let $S_i$ be the simple left module over $\la$ (resp., $\la_L$) that corresponds to $e_i$, i.e., $S_i$ is the $1$-dimensional module $Ke_i$ which is annihilated by the idempotents $e_j$ different from $e_i$ and by the ideal $\A$  (resp., by the canonical image of $\A$ in $\la_L$). The Jacobson radical of any $K$-algebra $\Delta$ we denote by $J(\Delta)$ or just $J$. Observe that $J(\la_L)$ is the canonical image of $\A$, whereas $J(\la)$ is the ideal generated by the paths $e_i \rightsquigarrow e_j$ between vertices $e_i, e_j$ not connected by  a path $e_j \rightsquigarrow e_i$ (see, e.g., \cite[Exercise (4), p.5]{CB.reps.quivers}). For any object $M$ of $\adm$, a \emph{top element of} $M$ is an element $x \in M \setminus \A M$ which satisfies $x = e_i x$ for some $i \le n$ (we then say that $x$ is \emph{normed by $e_i$}). A \emph{full set of top elements of $M$} is a set of top elements which generates $M$ and is linearly independent modulo $\A M$.

Given a dimension vector $\bd = (d_1, \dots, d_n)$  of $Q$, we let $\rep_{\bd}(\la)$, resp. $\rep_{\bd}(\la_L)$, be the classical affine variety parametrizing the (left) $\la$-modules (resp., $\la_L$-modules) with dimension vector $\bd$. Namely,
$$\rep_{\bd}(\la) \ = \ \prod_{\alpha \in Q_1} \Hom_K(K^{d_{\text{start}(\alpha)}},  K^{d_{\text{end}(\alpha)}} )$$
where $Q_1$ is the set of arrows of $Q$,
and
\begin{multline*}
\rep_{\bd}(\la_L) \ = \\
 \bigl\{(x_\alpha) \in \rep_\bd(\la) \bigm| x_{\alpha_L}\circ x_{\alpha_{L-1}} \circ \cdots \circ x_{\alpha_1} = 0 \ \text{for all paths}\ \alpha_L \cdots \alpha_1\ \text{in}\ Q \bigr\}.
\end{multline*}
These affine varieties come equipped with representation maps, i.e., with surjections sending any point $x$ of some $\rep_\bd(\la)$ to the isomorphism class of a module $M_x$ with dimension vector $\bd$; analogously for $\rep_\bd(\la_L)$.  The fibers of these maps are the orbits under the natural conjugation action of $\GL(\bd) := \prod_{1 \le i \le n} \GL_{d_i}$ on the displayed varieties. Finally, given any subset $U$ of one of these parametrizing varieties, we say that a module $M$ is ``in" $U$ if $M \cong M_x$ for some $x \in U$.  

\subsec*{Graphing. An informal introduction} 
We profile the structure of an admissible $\la$-module (or class of $\la$-modules) $M$ by (layered and labeled) graphs, that is, by finite undirected graphs in which each vertex represents a nonzero element of $M$, which is normed by some primitive idempotent $e_i$ and is correspondingly labeled by $i$. The set of vertices is partitioned into subsets $\L_0, \L_1, \dots, \L_k$, where each $\L_l$ consists of a full set of top elements of $\A^l M$ and $k$ is minimal with $\A^{k+1} M = 0$; we refer to the subsets $\L_j$ as the \emph{layers} of the graph and line up their elements in horizontal rows (the layers), where layer $0$ is displayed in the top row.  
Given that the label of any vertex coincides with that of the norming idempotent, the graph may contain arbitrarily high repetitions of any of the numbers $1, 2, \dots, n$ indexing the vertices of $Q$.  As for the edges: each is labeled by an arrow of $Q$; if $\alpha: e_i \rightarrow e_j$ is an arrow of $Q$, an edge labeled $\alpha$ from a vertex $x = e_i x$ in some layer $\L_l$ to a vertex $y = e_j y$ in some lower layer conveys the information $\alpha x \in K^*y$. 

Thus, in effect, these layered graphs are actually directed graphs, directed from the top down in our visual presentation, to allow for easy grasp of the simple composition factors of the consecutive quotients $\A^l M/ \A^{l+1} M$; arrow heads are omitted to distinguish graphs of modules from quivers of algebras.

For example, a graph of the form 
$$\xymatrixrowsep{1.7pc} \xymatrixcolsep{1pc}
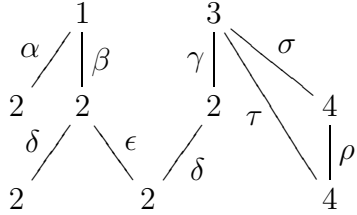
\xymatrix{
&1 \edge[dl]_{\alpha} \edge[d]^{\beta} &&3 \edge[d]_{\gamma} \edge[ddrr]_{\tau} \edge[drr]^{\sigma}  \\
2 &2 \edge[dl]_{\delta} \edge[dr]^{\epsilon} &&2 \edge[dl]^{\delta} &&4 \edge[d]^{\rho}  \\
2 &&2 &&&4 
}$$
profiles a (one-parameter) family of modules $M$ in which $M/\A M \cong S_1\oplus S_3$, \  $\A M/\A^2M \cong S_2^3 \oplus S_4$, and $\A^2 M \cong S_2^2 \oplus S_4$. Further, the graph indicates that  $M$ has a full set of top elements $z_1 = e_1 z_1$ and $z_2 = e_3 z_2$ such that 
\begin{align*}
\A M &= \la \alpha z_1 + \la \beta z_1 + \la \gamma z_2 + \la \sigma z_2\ \ \ \text{and}  \\
\A^2M &= \la \delta \beta z_1 + \la \epsilon \beta z_1 + \la \tau z_2 = \la \delta \beta z_1 + \la \epsilon \beta z_1 + \la \rho \sigma z_2  \\
&=  \la \delta \beta z_1 + \la \delta \gamma z_2 + \la \tau z_2 =  \la \delta \beta z_1 + \la \delta \gamma z_2 + \la \rho \sigma z_2  \, .
\end{align*}
The terminus of a downward path $p$ from a vertex representing a top element $z$ of $M$ stands for the element $pz$ of $M$ up to a nonzero scalar factor; e.g., the graph conveys that $\tau z_2$ and $\rho \sigma z_2$ are scalar multiples of each other, as are $\epsilon \beta z_1$ and $\delta \gamma z_2$ (the latter scalar may be normalized to $1$ by passing to a suitable scalar multiple of one of the top vertices $z_1$, $z_2$).  That no edge labeled $\epsilon$ connects to the edge labeled $\gamma$ indicates that $\epsilon \gamma z_2 = 0$. In particular, this means that the elements of $M$ represented by the vertices of the graph form a $K$-basis of $M$.

Not every admissible $\la$-module $M$ has a graph as described above.  To fill the gap, ``pooling" is used to signal linear dependencies of more than two elements $p_i z_i$ in $M$.  Suppose $r \ge 3$. Let $z_1, \dots, z_r$ be top vertices of the graph of $M$ and $p_1, \dots, p_r$ paths in $Q$ ending in a vertex $e_i$ such that the set $\{p_1 z_1, \dots, p_r z_r\}$ is minimally linearly dependent in the sense that this set is linearly dependent while every $(r-1)$-element subset is linearly independent.  Then the requirement that the vertices should form a $K$-basis for $M$ is modified as follows:  the minimal linear dependence is indicated by a dotted line encircling the vertices labeled $i$ that represent $p_1 z_1, \dots, p_r z_r$.

For instance, the elements $\alpha z_1$, $\beta z_1$, $\delta \beta z_1$ are linearly independent in any module of the class profiled by the graph above, whereas adding a pool to the graph as below conveys that the $K$-subspace spanned by $\alpha z_1$, $\beta z_1$, $\delta \beta z_1$ has dimension $2$:
$$\xymatrixrowsep{1.7pc} \xymatrixcolsep{1pc}
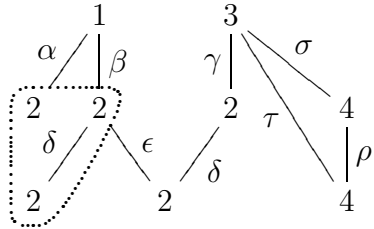
\xymatrix{
&1 \edge[dl]_{\alpha} \edge[d]^{\beta} &&3 \edge[d]_{\gamma} \edge[ddrr]_{\tau} \edge[drr]^{\sigma}  \\
2 \save[0,0]+(0,3);[0,0]+(0,3) **\crv{~*=<2.5pt>{.} [0,0]+(3,3) &[0,1]+(0,3) &[0,1]+(3,3) &[0,1]+(3,0) &[1,0]+(3,-3) &[1,0]+(0,-3) &[1,0]+(-3,-3) &[1,0]+(-3,0) &[0,0]+(-3,0) &[0,0]+(-3,3)} \restore &2 \edge[dl]_{\delta} \edge[dr]^{\epsilon} &&2 \edge[dl]^{\delta} &&4 \edge[d]^{\rho}  \\
2 &&2 &&&4
}$$

When $Q$ does not have multiple equi-directed arrows between any two vertices, the edge labels in the graph of a $\la$-module may be omitted without loss of information.  In this vein, the labels $\gamma$, $\rho$, $\sigma$, $\tau$ could be omitted in the example above.

\sectionnew{\bf Irreducible components of the parametrizing varieties and their generic modules}

\noindent  Let $\bd$ be a dimension vector of $Q$. In this section we focus on the subvarieties $\rep_{\bd}(\la_L)$ of $\rep_\bd(\la)$ for $L \in \NN$.  

\begin{observation}  \label{obs1} 
$\rep_{\bd}(\la)$ is irreducible for any dimension vector $\bd$, whereas the $\rep_{\bd} (\la_L)$ for $L \ge 2$ are reducible in general.
\medskip

Indeed, the first claim is immediate from the fact that $\rep_\bd(\la) \cong \AA^N$, where $N = {\sum_{\alpha  \in Q_1} d_{\text{start}(\alpha)}\cdot d_{\text{end}(\alpha)}}$. For the second consider, e.g., the quiver 
$\xymatrixrowsep{3pc} \vcenter{\xymatrix{ 1 \ar@/^/[r] &2 \ar@/^/[l] }}$
and $\bd = (1,1)$; for any $L \ge 2$, the variety $\rep_\bd(\la_L)$ then has two irreducible components, one consisting of the uniserial module with top $S_1$ next to the semisimple module $S_1 \oplus S_2$, the other consisting of the uniserial module with top $S_2$ next to $S_1 \oplus S_2$.  In fact, based on the oriented cycle with $L\ge 2$ vertices and $L$ arrows, and $\bd = (1, 1, \dots, 1)$, one analogously obtains an example where $\rep_\bd(\la_L)$ has precisely $L$ irreducible components.
\end{observation}

\noindent 
Thus, in a first step, the plan to explore generic properties of the representations pa\-ra\-met\-rized by $\rep_\bd(\la_L)$ calls for a useful listing of the irreducible components.

\subsec*{Outline of Section 2} In Subsection 2.A, we  introduce and discuss  concepts that are relevant to a classification of the components of $\rep_\bd(\la_L)$; in Subsection 2.B, we  state the main results; Subsection 2.C  targets generic modules in the components and their generic minimal projective presentations; Subsection 2.D  illustrates the theory and addresses its algorithmic aspect.

\subsec*{2.A. Preparing for the main results } 

The problem of understanding the irreducible components of $\rep_{\bd}(\Delta)$ for a finite dimensional algebra $\Delta$ has been tackled intensely.  We do not include an overview of the status quo, but  restrict our attention to the case where $\Delta$ is a (truncated) path algebra. 

We begin with a few general comments. Suppose we already have an irreducible subset $U$ of $\rep_\bd(\Delta)$.  Particularly useful tools for deciding whether the closure $\overline{U}$ of $U$ is a component come in the form of collections of upper or lower semicontinuous maps $f$ from $\rep_{\bd}(\Delta)$ to partially ordered sets. In the case of upper semicontinuity, attainment of a minimal value of $f$ on $U$ indicates that $\overline{U}$ is an irreducible component of $\rep_{\bd}(\Delta)$, whereas attainment of a maximal value on $U$ leads to the same conclusion if $f$ is lower semicontinuous.  Examples of lower semicontinuous maps on $\rep_\bd(\Delta)$ are path ranks, i.e., maps of the form $x \mapsto \dim p M_x$ where $p$ is a path in the quiver of $\Delta$; or the map that assigns to each $x \in \rep_\bd(\Delta)$ the dimension of the $\GL(\bd)$-orbit of $x$.  Instances of upper semicontinuous maps, our tool of preference in tackling $\Delta = \la_L$, are $\dim \Hom_\Delta(M_x, N)$, $\dim \Hom_\Delta(N, M_x)$, $\dim \Ext^i_\Delta(M_x, N)$, and $\dim \Ext^i_\Delta(N, M_x)$ for $N \in \deltamod$.  However, the most effective towards the detection of components in our present setting are upper semicontinuous maps based on filtrations of the modules $M_x$, which we  describe next.  

Now assume that $\Delta$ has Loewy length $\le L$.
Our reference filtration of a module $M \in \deltamod$ is its \emph{radical filtration}, namely the chain
$$M \supseteq JM \supseteq J^2 M  \supseteq \cdots \supseteq J^L M = 0.$$ 
The family of successive quotients $\bigl(\SS_l(M)\bigr)_{0 \le l < L}$, where $\SS_l(M) = J^l M / J^{l+1} M$, is referred to as the \emph{radical layering of $M$}; clearly $\SS(M) = (\SS_0(M), \SS_1(M), \dots, \SS_{L-1}(M))$  is then an $L$-term sequence of semisimple modules whose dimension vectors add up to 
$$\udim \SS(M): = \sum_{0 \le l \le L-1} \udim \SS_l(M) = \udim M.$$  

Any sequence $\SS = (\SS_0, \SS_1, \dots, \SS_{L-1})$ of semisimple modules with $\udim \SS = \bd$ is called a \emph{semisimple sequence with dimension vector $\bd$}.  We  identify isomorphic semisimple modules, whence the set $\fancyS(\bd)$ of all semisimple sequences with dimension vector $\bd$ is finite. 

Another natural example of a semisimple sequence comes from the \emph{socle filtration} 
$$\soc_0 (M) := 0  \subseteq \soc_1(M) := \soc(M) \subseteq \soc_2(M) \subseteq \cdots \subseteq \soc_L(M) = M,$$
which gives rise to the \emph{socle layering} $\SS^*(M) := \bigl( \SS^*_l(M) \bigr)_{0\le l<L}$ where
$$
\SS^*_l(M) = \soc_{l+1}(M)/ \soc_l(M) = \soc \bigl( M/\soc_l(M) \bigr).
$$
On equipping the set $\fancyS(\bd)$ with the following dominance partial order, the assignments
$$x \longmapsto \SS(M_x) \qquad\text{and} \qquad x \longmapsto \SS^*(M_x)$$
yield upper semicontinuous maps $\rep_{\bd}(\Delta) \longrightarrow \fancyS(\bd)$.
The \emph{dominance order on $\fancyS(\bd)$} is as follows:  
$$\SS \le \widetilde{\SS} \quad \iff \quad \udim \bigoplus_{0 \le l \le m} \SS_l \le \udim \bigoplus_{0 \le l \le m} \widetilde{\SS}_l \quad \forall\ m \in \{0, \dots, L-1\},$$
where the partial order on dimension vectors is the componentwise order.

A prominent role is played by the semisimple sequences which are realizable as radical layerings:

\begin{definition}  \label{realizable}
A semisimple sequence $\SS$ is called \emph{realizable} (with respect to $\deltamod$) if $\SS = \SS(M)$ for some $\Delta$-module $M$. 
\end{definition}

\begin{observation}  \label{obs2}
Both of the maps $\rep_{\bd}(\Delta) \longrightarrow \fancyS(\bd)$ given by $x \mapsto \SS(M_x)$ and $x \mapsto \SS^*(M_x)$ are upper semicontinuous \cite[Observation 2.10]{irredcomp.local}. Hence so is the pairing $\Theta: \rep_{\bd}(\Delta) \longrightarrow \fancyS(\bd)^2$ defined by $\Theta(x) = \bigl(\SS(M_x), \SS^*(M_x) \bigr)$.
\end{observation} 

\noindent 
In numerous special cases, the map $\Theta$ detects all irreducible components of the varieties $\rep_\bd(\Delta)$, meaning that for any irreducible component $\C$ of $\rep_\bd(\Delta)$, there is some $x \in \C$ such that $\Theta(x)$ is minimal in $\Im(\Theta)$.  This holds, for instance, in the truncated local case, i.e., when $\Delta = \la_L$ and $Q$ consists of a single vertex with two or more loops \cite[Main Theorem]{irredcomp.local}. However, in general this map still has blind spots even for $\Delta = \la_L$ (see, e.g., \cite[Example 4.8]{irredcomp.local}). 

On the other hand, there does exist an upper semicontinuous map which, in case $\Delta = \la_L$, sees \emph{all} irreducible components of the varieties $\rep_\bd(\Delta)$ for arbitrary quivers $Q$. It turns out that the key to detecting all components in this case is provided by module filtrations which are governed by a realizable semisimple sequence in the following sense.  

\begin{definition}  \label{governed}
\cite[Definitions 3.5, 3.10]{closures} Let $\Delta = KQ/ I$ be any (finite dimensional) path algebra modulo relations, and $M \in \deltamod$ with dimension vector $\bd$. Moreover, suppose that $\SS = \bigl(\SS_l\bigr)_{0 \le l \le L-1}$ is a semisimple sequence of $\Delta$-modules with dimension vector $\bd$. 
\smallskip

$\bullet$ A filtration $M = M_0 \supseteq M_1 \supseteq M_2 \supseteq \cdots \supseteq M_L = 0$ is said to be \emph{governed by $\SS$} if $M_l/ M_{l+1} \cong \SS_l$ for $0 \le l \le L-1$.
\smallskip

$\bullet$ Now suppose that the semisimple sequence $\SS$ is realizable. We define $\Filt \SS $ to be the (closed) set of those points $x \in \rep_\bd(\Delta)$ which have the property that $M_x$ has at least one filtration governed by $\SS$.
\smallskip

$\bullet$ Let $\Gamma: \rep_{\bd}(\Delta) \longrightarrow \NN$ be the map sending any point $x \in \rep_{\bd}(\Delta)$ to the number of \emph{realizable} semisimple sequences that govern some filtration of $M_x$.  (Note that $\Gamma$ is well-defined: indeed, $\Gamma(x) \ge 1$, because the radical layering $\SS(M_x)$, a realizable semisimple sequence, governs the radical filtration of $M_x$.)  
\end{definition}

\begin{theorem}  \label{obs3}
{\rm\cite[Corollary 3.11]{closures}}
For any finite dimensional path algebra modulo relations $\Delta = KQ/ I$ and any dimension vector $\bd$ of $Q$, the map $\Gamma:  \rep_{\bd}(\Delta) \longrightarrow \NN$ is upper semicontinuous.
\end{theorem}
\smallskip

\subsec*{2.B. Classifying the irreducible components of  $\rep_\bd(\la_L)$ and tracking their evolution for growing $L$.} 
We are now in a position to state the main results of the section.

\begin{theorem}
{\rm\cite[Observations 2.1, 2.4, Theorem 5.3, Proposition 2.2]{BHT}}  For any quiver $Q$, any integer $L \ge 2$, and any realizable semisimple sequence $\SS$ with dimension vector $\bd$, the following locally closed subvariety of $\rep_\bd(\la_L)$ is irreducible:
$$\rep \SS: = \{x \in \rep_\bd(\la_L) \ \mid \ M_x \ \text{has radical layering} \ \SS\}.$$
Moreover, all irreducible components of $\rep_\bd(\la_L)$ are among  the Zariski closures of the subvarieties $\rep \SS$, where $\SS$ traces the realizable semisimple sequences with dimension vector $\bd$.
\end{theorem}

\noindent 
Determining the irreducible components of $\rep_\bd(\la_L)$ therefore boils down to identifying those realizable semisimple sequences $\SS$ which give rise to maximal irreducible subsets $\overline{\rep \SS}$ of $\rep_\bd(\la_L)$. It turns out that the upper semicontinuous map $\Gamma$ introduced in the preceding subsection sifts out all of them. 
\bigskip

\begin{theorem}  \label{repSbar}
{\rm\cite[Theorems 4.3, 4.5, \S5B]{closures}, \cite[Criterion 3.2]{irredcomp.local}} Let $\SS$ be a realizable semisimple sequence with dimension vector $\bd$.  Then 
\smallskip

\noindent {\bf (1)} The closure $\overline{\rep \SS}$ is an irreducible component of $\rep_\bd(\la_L)$ if and only if $1 \in \Gamma(\rep \SS)$.  Thus, $\SS$ gives rise to an irreducible component of $\rep_\bd(\la_L)$ precisely when there exists a module $M$ in $\rep \SS$ such that the only filtration of $M$ governed by a realizable semisimple sequence is its radical filtration.
\medskip

\noindent {\bf (2)} $\overline{\rep \SS} = \Filt \SS$ = the subvariety of $\rep_\bd(\la_L)$ consisting of the points $x$ such that $M_x$ has a filtration governed by $\SS$. 
\medskip

\noindent {\bf (3)} The irreducible components of $\rep_\bd(\la_L)$ can be algorithmically determined from $Q$ and $L$. Moreover, given irreducible components $\C_i$ of $\rep_\bd(\la_L)$ and $\D_j$ of $\rep_\bd(\la_{L'})$ for $L < L'$, the question whether $\C_i$ is contained in $\D_j$ can be algorithmically decided.
\end{theorem}

\noindent 
In general the available algorithms as in part (3) of the preceding theorem are labor-intensive.  However, in case $Q$ is the quiver consisting of a single vertex and $2$ or more loops, the irreducible components of the varieties $\rep_\bd(\la_L)$ can be listed without any computational  effort; this will become clear from the last of our main results.   In this case, it is also less laborious to establish the full inclusion hierarchy of the components of the $\rep_\bd(\la_L)$ as $L$ grows; cf.~Example \ref{Qa1b}.   

\begin{theorem}  \label{mainthm.local}
{\rm\cite[Main Theorem]{irredcomp.local}}
Suppose $Q$ is the quiver
$$\xymatrixrowsep{0.75pc}\xymatrixcolsep{0.3pc}
\xymatrix{
 & \ar@{{}{*}}@/^1pc/[dd] \\
1 \ar@'{@+{[0,0]+(-6,-6)} @+{[0,0]+(-15,0)} @+{[0,0]+(-6,6)}}^(0.7){\alpha_1}
\ar@'{@+{[0,0]+(-6,6)} @+{[0,0]+(0,15)} @+{[0,0]+(6,6)}}^(0.7){\alpha_2}
\ar@'{@+{[0,0]+(6,-6)} @+{[0,0]+(0,-15)} @+{[0,0]+(-6,-6)}}^(0.3){\alpha_r}  \\
 &
}$$
where $r \ge 2$, and $d \in \NN$.  
\smallskip

{\bf(I)} If $L < d$, the following conditions are equivalent for any realizable semisimple sequence $\SS$ with dimension $d$:
\smallskip

\quad {\bf (1)} $\overline{\rep \SS}$ is an irreducible component of $\rep_d(\la_L)$.
\smallskip

\quad {\bf (2)} $\dim \SS_l \le r \cdot \dim \SS_{l-1}$ and $\dim \SS_{l-1} \le r \cdot \dim \SS_l$ for $1 \le l \le L-1$.
\smallskip

\quad {\bf (3)} The semisimple sequence $\SS^* = (\SS_{L-1}, \dots, \SS_0)$ is the generic socle layering of the 

$\quad$ \ \ \ \ \ modules in $\rep \SS$.
\vskip.03truein

\noindent As is the case for an arbitrary quiver $Q$, all irreducible components of $\rep_d(\la_L)$ are among the $\overline{\rep \SS}$.
\smallskip

{\bf(II)} If $L \ge d$, the variety $\rep_d(\la_L)$ is irreducible and, generically, its modules are uniserial.
\end{theorem}

\subsec*{2.C. Generic modules for the irreducible components of the varieties $\rep_\bd(\la_L)$}
Here ``generic" is used in the geometric sense of being in suitably general position within a component, rather than being indecomposable and infinite dimensional with finite endolength.  More precisely, we consider $\SS$-generic modules, defined as follows.  Assume $K$ has infinite transcendence degree over an algebraic closure $K_0$ of its prime subfield.  A $\la_L$-module $G$ in some $\rep\SS$ is \emph{$\SS$-generic} if $G$ satisfies all categorically defined module properties which (i) hold for all modules in some dense open subset of $\SS$ and (ii) are preserved by the Morita self-equivalences of $\lalmod$ induced by the $K_0$-automorphisms of $K$.

For the following results we summarily refer to \cite[Sections 4,5]{BHT} and to the slimmer presentation in \cite[Section 2]{irredcomp.acyclic}.  We simplify more general statements by only formulating them for truncated path algebras $\la_L$.

\begin{theorem}  
{\rm\cite[Theorems 4.3, 5.12]{BHT}}
Suppose that the base field $K$ has infinite transcendence degree over its prime field, and let $\SS$ be any realizable semisimple sequence in $\la_L\mbox{\rm-mod}$.  Then there exists an $\SS$-generic module $G$ in $\rep \SS$, and $G$ has an algorithmically obtainable minimal projective presentation.
\end{theorem}

\begin{observation} (Loosely phrased; see \cite[Section 2.B]{irredcomp.acyclic})  If $K$ fails to have infinite transcendence degree over its prime field, one passes to $\widehat{\la} = \widehat{K} Q$, where $\widehat{K}$ is the algebraic closure of the rational function field $K(X_n \mid n \in \NN)$. In the context of our present investigation, this is harmless because the exact and faithful functor $\widehat{K} \otimes_K - : \lamod \longrightarrow \widehat{\la}\text{-mod}$ preserves and reflects dimension vectors of modules and semisimple sequences as well as the generic module properties in which we are interested.
\end{observation}

\noindent 
Thus, in deciding whether $\rep \widetilde{\SS} \subseteq \overline{\rep \SS}$, we only need to check whether a generic module for $\rep \widetilde{\SS}$ has a filtration governed by $\SS$. 

\begin{example}  \label{Q1a2bb1}
Let $Q$ be the quiver below, $\la = \CC Q$, and $L = 3$, $\bd = (2,2)$. 
$$\xymatrixcolsep{3pc}
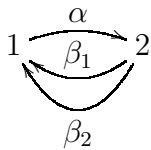
\xymatrix{
1 \ar@/^/[r]^{\alpha} &2 \ar@/^1pc/[l]_{\beta_1} \ar@/^2pc/[l]^{\beta_2}
}$$ 

\noindent Consider the realizable semisimple sequences
$$
\SS^{(1)} = (S_1 \oplus S_2, S_1 \oplus S_2, 0)), \quad\SS^{(2)} = (S_1 \oplus S_2, S_1, S_2), \quad \SS^{(3)} = (S_2, S_1^2, S_2).
$$
The corresponding varieties $\rep \SS^{(j)}$ have generic modules $G_j$ with the following graphs:
$$\xymatrixrowsep{1.5pc}\xymatrixcolsep{0.6pc}
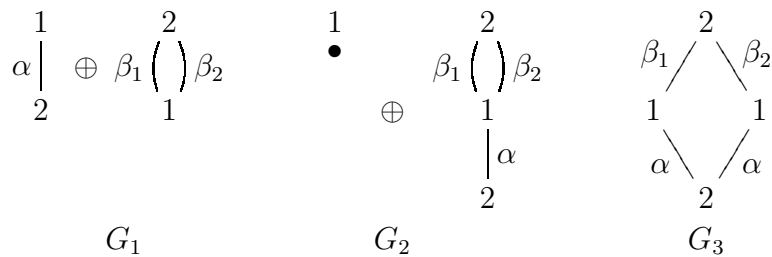
\xymatrix{
1 \edge[d]_{\alpha} &\ar@{}[d]|{\oplus} &&2 \edge@/_/[d]_{\beta_1} \edge@/^/[d]^{\beta_2} & &&&1 \drbl &&&2 \edge@/_/[d]_{\beta_1} \edge@/^/[d]^{\beta_2} & &&& &2 \edge[dl]_{\beta_1} \edge[dr]^{\beta_2}  \\
2 &&&1 &&& &&\oplus &&1 \edge[d]^{\alpha} & &&&1 \edge[dr]_{\alpha} &&1 \edge[dl]^{\alpha}  \\
&&\dropdown{G_1} &&& &&  &\dropdown{G_2}  &&2 &&& &&2 \dropdown{G_3}
}$$
Then $\rep \SS^{(i)} \subseteq \overline{\rep \SS^{(3)}}$ for $i = 1,2$.  Indeed, $G_1$ has a chain of submodules
$$U_0 = G_1 \ \supset \ U_1 = {\vcenter{\xymatrixrowsep{1.3pc}\xymatrixcolsep{0.6pc}\xymatrix{1 \edge[d] \\ 2}}} \oplus S_1 \ \supset \ U_2 = S_2 \ \supset \ U_3 = 0,$$
and it is equally easily checked that also $G_2$ has a filtration governed by $\SS^{(3)}$.
\end{example}

\subsec*{2.D. Further examples}

In carrying out computations towards a full list of irreducible components of the varieties $\rep_\bd(\la_L)$ for specific quivers, it is advantageous to secure as many components as possible by way of the upper semicontinuous map $\Theta$ of subsection 2.A, or by combining $\Theta$  with the lower semicontinuous maps $\text{nullity}_p: \rep_\bd(\la_L) \longrightarrow \ZZ_{\ge0}$ for the paths $p$ of $Q$, which send $x$ to the nullity of $M_x \longrightarrow M_x$, $y \mapsto p y$, respectively.  Testing for existence of $\SS$-governed filtrations carries the highest labor cost.

\begin{example}  \label{Q1a2bb1.more} We return to the concluding example of subsection 2.C, i.e., 
we let $\la = \CC Q$ and $\bd = (2,2)$, where $Q$ is the quiver of Example \ref{Q1a2bb1}. We first state the outcome regarding the components of $\rep_\bd(\la_L)$ for $L = 3,4$, and then justify some samples. Since $|\bd| = 4$, we know that the sequence $\rep_\bd(\la_L)$, $L \in \NN$, becomes stationary for $L \ge 4$.

When $L=3$, there are $4$ irreducible components, with generic modules having graphs of the form
$$\xymatrixrowsep{1.5pc}\xymatrixcolsep{0.5pc}
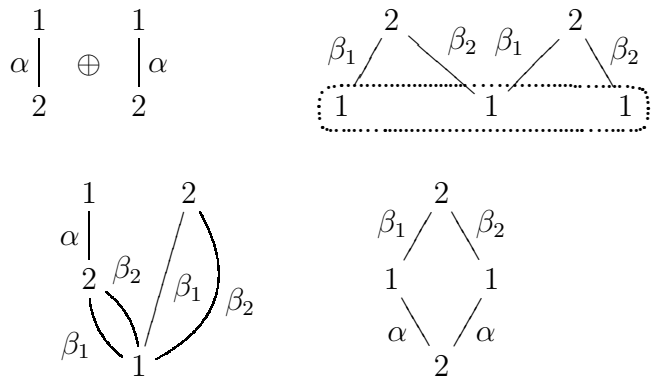
\xymatrix{
1 \edge[d]_{\alpha} &\ar@{}[d]|{\oplus} &1 \edge[d]^{\alpha} &&&&& &2 \edge[dl]_{\beta_1} \edge[drr]^{\beta_2} &&&&2 \edge[dll]_{\beta_1} \edge[dr]^{\beta_2}  \\
2 &&2 &&&&&1 \horizpool{6} &&&1 &&&1  \\  
&1 \edge[d]_{\alpha} &&2 \edge[ddl]^{\beta_1} \edge@/^1.7pc/[ddl]^{\beta_2} &&&&& &2 \edge[dl]_{\beta_1} \edge[dr]^{\beta_2}  \\
&2 \edge@/^/[dr]^(0.2){\beta_2} \edge@/_/[dr]_{\beta_1} && &&&&&1 \edge[dr]_{\alpha} &&1 \edge[dl]^{\alpha}  \\
&&1 & &&&&& &2
}$$

When $L\ge4$, there are just $2$ irreducible components,
$$
\overline{\rep(S_1,S_2,S_1,S_2)} \qquad \text{and} \qquad \overline{\rep(S_2,S_1,S_2,S_1)},
$$
both containing infinitely many orbits of uniserial modules with dimension $4$.

As far as containments of the components for $L=3$ in  components for $L=4$ are concerned:  The component in the northeast corner of the display for $L=3$ is contained in the closure of $\rep(S_1, S_2, S_1, S_2)$ in $\rep_\bd(\la_4)$ and in no others, while the remaining components for $L=3$ are  contained in the closure of $\rep(S_2,S_1,S_2,S_1)$ only.

\emph{Partial justifications}: To verify the list of irreducible components for $L = 3$, observe that, among the realizable semisimple sequences $\SS^{(1)}$, $\SS^{(2)}$, $\SS^{(3)}$ listed in Example \ref{Q1a2bb1}, only $\SS^{(3)}$ is eligible to be the generic radical layering of an irreducible component of $\rep_\bd(\la_3)$. There are precisely $3$ more realizable semisimple sequences in this case, namely $\SS^{(4)} = (S_1^2, S_2^2)$, $\SS^{(5)} = (S_2^2, S_1^2)$, and $\SS^{(6)} = (S_1 \oplus S_2, S_2, S_1)$ with generic modules as shown above. Clearly, none of their generic socle layerings is comparable with that of $G_3$, whence $\rep \SS^{(3)} = \rep (S_2, S_1^2, S_2)$ is indeed an irreducible component. Now focus on $\SS^{(6)}$. The map $\Theta$ rules out the possibility that $\rep \SS^{(6)}$ be contained in the closure of $\rep \SS^{(j)}$ for $j = 4$ or $5$. To exclude that $\rep \SS^{(6)}$ be contained in the closure of $\rep \SS^{(3)}$, we check that the displayed generic module $G_6$ for $\rep \SS^{(6)}$ does not have a filtration governed by $\SS^{(3)}$:  Clearly $M_0 = G_6$ has a submodule $M_1$ with the property that $M_0/M_1 \cong S_2$.  However, no eligible choice of $M_1$ has two copies of $S_1$ in its top, which confirms the claim.    \qed
\end{example}

\noindent 
We next illustrate Theorem \ref{mainthm.local} for $r = 2$.  

\begin{example}  \label{Qa1b}
Let $Q$ be the quiver \ $\vphantom{\biggl(}\xymatrix{ \bullet \ar@(ul,dl)_{\alpha} \ar@(ur,dr)^{\beta} }$, \ and $d=5$.
 
The generic radical layerings of the irreducible components of $\rep_5(\la_L)$, for $L = 2,\dots,5$, are as follows.  In the diagram for each component, the number of bullets in a given layer indicates the dimension of that radical layer, i.e., the multiplicity of the unique simple module in that layer.
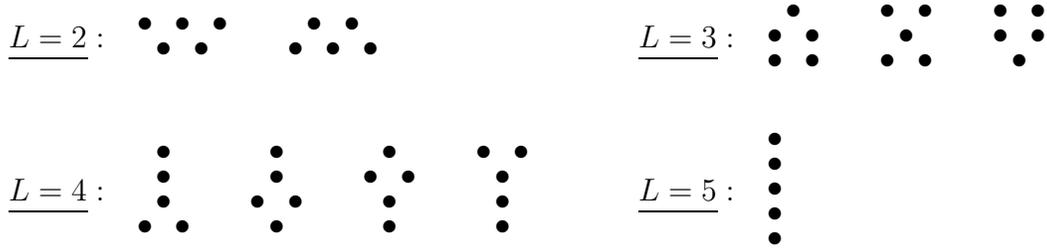
\begin{align*}
\underline{\vphantom{p}{L=2}} : &\quad\xymatrixrowsep{0.2pc}\xymatrixcolsep{0.01pc} 
\vcenter{\xymatrix{
\pbl &&\pbl &&\pbl &&&& &\pbl &&\pbl  \\
&\pbl &&\pbl & &&&&\pbl &&\pbl &&\pbl
}}
&\quad \underline{\vphantom{p}{L = 3}} : &\quad\xymatrixrowsep{0.2pc}\xymatrixcolsep{0.01pc} 
\vcenter{\xymatrix{
&\pbl & &&&&\pbl &&\pbl &&&&\pbl &&\pbl  \\
\pbl &&\pbl &&&& &\pbl & &&&&\pbl &&\pbl  \\
\pbl &&\pbl &&&&\pbl &&\pbl &&&& &\pbl
}}  \\  \\
\underline{\vphantom{p}{L=4}} : &\quad\xymatrixrowsep{0.2pc}\xymatrixcolsep{0.01pc} 
\vcenter{\xymatrix{
&\pbl & &&&& &\pbl & &&&& &\pbl & &&&&\pbl &&\pbl  \\
&\pbl & &&&& &\pbl & &&&&\pbl &&\pbl &&&& &\pbl  \\
&\pbl & &&&&\pbl &&\pbl &&&& &\pbl & &&&& &\pbl  \\
\pbl &&\pbl &&&& &\pbl & &&&& &\pbl & &&&& &\pbl
}}
&\quad \underline{\vphantom{p}{L=5}} : &\quad\xymatrixrowsep{0.2pc}\xymatrixcolsep{0.01pc} 
\vcenter{\xymatrix{
\pbl \\  \pbl \\  \pbl \\  \pbl \\  \pbl
}}
\end{align*}

\noindent Let us label the corresponding irreducible components as follows: $\C_1$, $\C_2$ (for $L=2$); $\D_1$, $\D_2$, $\D_3$ (for $L=3$); $\E_1$, $\E_2$, $\E_3$, $\E_4$ (for $L=4$); and $\F$ (for $L=5$).  The inclusion hierarchy of the components of the $\rep_\bd(\la_L)$ for increasing $L$ is shown below.
\smallskip

$$\xymatrixrowsep{2.5pc}\xymatrixcolsep{1pc}
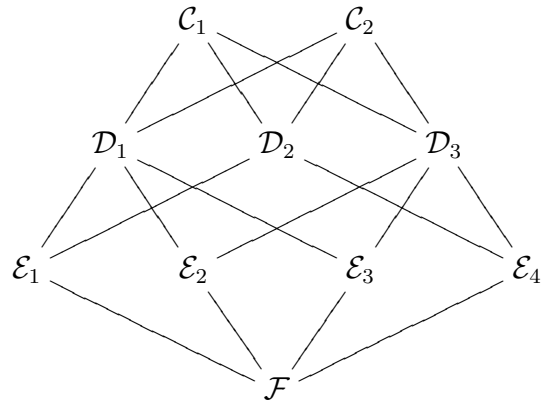
\xymatrix{
&&\C_1 \edge[dl] \edge[dr] \edge[drrr] &&\C_2 \edge[dlll] \edge[dl] \edge[dr]  \\
&\D_1 \edge[dl] \edge[dr] \edge[drrr] &&\D_2 \edge[dlll] \edge[drrr] &&\D_3 \edge[dlll] \edge[dl] \edge[dr]  \\
\E_1 \edge[drrr] &&\E_2 \edge[dr] &&\E_3 \edge[dl] &&\E_4 \edge[dlll]  \\
&&&\F
}$$
\smallskip

\noindent We give two samples of the arguments backing this containment graph, setting $J = J(\la_3)$ and denoting the unique simple left $\la_3$-module by $S$. 

To see that $\D_2$ is not contained in $\E_2$, start by observing that the generic radical layering $\widetilde{\SS} = (S^2, S, S^2, 0)$ of $\D_2$ is larger than or equal to the generic radical layering $\SS = (S, S, S^2, S)$ of $\E_2$, and the analogous inequality holds for the generic socle layerings of $\D_2$ and $\E_2$; hence the test map $\Theta$ does not rule out the possibility that $\D_2$ be contained in $\E_2$. Letting $M_0$ be a generic module for $\D_2$, we  show that $M_0$ does not have a filtration governed by $\SS$. For that purpose, we display a generic graph of $M_0$. 
$$\xymatrixrowsep{1.5pc}\xymatrixcolsep{0.5pc}
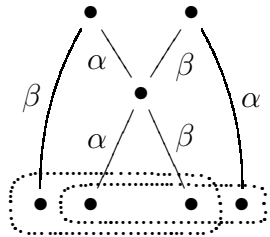
\xymatrix{
&\bullet \edge@/_/[ddl]_{\beta} \edge[dr]_{\alpha} &&\bullet \edge[dl]^{\beta} \edge@/^/[ddr]^{\alpha}  \\
&&\bullet \edge[dl]_{\alpha} \edge[dr]^{\beta}  \\
\bullet \save[0,0]+(0,4);[0,0]+(0,4) **\crv{~*=<2.5pt>{.} [0,3]+(0,4) &[0,3]+(3,4) &[0,3]+(4,0) &[0,3]+(4,-4) &[0,3]+(0,-4) &[0,0]+(0,-4) &[0,0]+(-4,-4) &[0,0]+(-4,0) &[0,0]+(-4,4)}\restore
 &\bullet \save[0,0]+(0,2.5);[0,0]+(0,2.5) **\crv{~*=<2.5pt>{.} [0,3]+(0,2.5) &[0,3]+(4,2.5) &[0,3]+(2.5,0) &[0,3]+(4,-2.5) &[0,3]+(0,-2.5) &[0,0]+(0,-2.5) &[0,0]+(-4,-2.5) &[0,0]+(-4,0) &[0,0]+(-4,2.5)}\restore &&\bullet &\bullet &\vphantom{\biggl(}
}$$
Clearly, any submodule $M_1 \subset M_0$ with $M_0/ M_1 \cong S$ properly contains $J M_0$ and has top $S$. Consequently, the only submodule $M_2 \subset M_1$ with $M_1/ M_2 \cong S$ is $M_2 = JM_0$, which is in turn a module with top $S$.  Therefore, $M_2$ fails to have a submodule $M_3$ such that $M_2/ M_3 \cong S^2$.  This confirms that $M_0$ does not have a filtration governed by $\SS$, whence part (2) of Theorem \ref{repSbar} tells us that, indeed, $\D_2 \not\subseteq \E_2$.

Next we verify that $\D_3 \subseteq \E_3$. First we note that the test map $\Theta$ does not rule this out. Now we let  $M_0$ be a generic module for $\D_3$. The module $M_0/ J^2 M_0$, generic for the semisimple sequence $(S^2, S^2)$ has a graph of the form 
$$\xymatrixrowsep{1.5pc}\xymatrixcolsep{1.5pc}
\xymatrix{
\bullet \edge@/_/[d]_{\alpha} \edge@/^/[d]^{\beta} &\edge@{}[d]|{\oplus} &\bullet \edge@/_/[d]_{\alpha} \edge@/^/[d]^{\beta}  \\
\bullet &&\bullet
}$$

\noindent Hence there is a submodule $M_1\subset M_0$ such that $M_0/M_1 \cong S$ and $M_1/ JM_1 \cong S^2$. Finally, $M_2: = J M_1$ has a filtration governed by $(S, S)$, which completes the argument.
\end{example} 

\sectionnew{\bf The strong tilts of the algebras $\la_L$}

\subsec*{3.A. Generalities}
Here a \emph{tilt} of a finite dimensional algebra $\Delta$ is short for ``an algebra obtained from $\Delta$ by tilting", i.e., an algebra $\widetilde{\Delta} = \End_\Delta(T)^{\op}$, where $_\Delta T$ is a tilting module in the sense of Miyashita.  This means that $\pdim_\Delta T$ is finite, $\Ext_\Delta^i(T,T) = 0$ for $i \ge 1$, and there exists an exact sequence $0 \rightarrow {_\Delta}\Delta \rightarrow T^{(0)} \rightarrow T^{(1)} \rightarrow \cdots \rightarrow T^{(m)} \rightarrow 0$ with all $T^{(j)}$ in $\add(T)$, the category of finite direct sums of direct summands of $T$. It is well known that $_{\Delta}T_{\widetilde{\Delta}}$ is then a balanced bimodule, which is a tilting module in both $\Delta\text{-mod}$ and $\text{mod-}\widetilde{\Delta}$.

Thus tilting modules are mongrels of projective generators and injective cogenerators. Reflecting this fact, the corresponding Hom and tensor functors yield generalizations of both Morita equivalences (see, e.g., \cite{handbook}) and Morita dualities (see \cite{Miy} and \cite{dualities}).  More precisely, tilting modules induce equivalences between suitable pairs of subcategories of  $\Delta\text{-mod}$  and  $\widetilde{\Delta}$-mod and dualities between suitable pairs of subcategories of $\Delta\text{-mod}$ and mod-$\widetilde{\Delta}$. The strong homological connection between the algebras $\Delta$ and $\widetilde{\Delta}$ is witnessed by the fact that these algebras are derived equivalent (see \cite{Hap} and \cite{CPS}). 

The tilting modules closest to projective generators are obviously the projective generators themselves. In particular, there is a unique basic one.  As for a dual, the tilting modules ``closest to injective cogenerators" were singled out and dubbed \emph{strong tilting modules} by Auslander and Reiten \cite{AuRe}.  These are tilting modules which are injective cogenerators within the category $\pinf(\Delta\text{-mod})$ of finitely generated $\Delta$-modules of finite projective dimension. However, $\lamod$ need not have a strong tilting module.

\begin{theorem}  \label{strong.tilting}
{\rm (See \cite[Section 6]{AuRe} for part (1) and \cite[\S 2.A]{dualities} for part (2).)}  Let $\Delta$ be a finite dimensional algebra.
\smallskip

\noindent {\bf (1)} $\Delta\mbox{\rm-mod}$ has a strong tilting module if and only if the category $\pinf(\Delta\mbox{\rm-mod})$ is contravariantly finite in $\Delta\mbox{\rm-mod}$, meaning that every finitely generated left $\Delta$-module $M$ has a right approximation in $\pinf(\Delta\mbox{\rm-mod})$, i.e., a map $\varphi \in \Hom_\Delta(\pinf(\Delta\mbox{\rm-mod}), M)$ with the property that all maps in $\Hom_\Delta(\pinf(\Delta\mbox{\rm-mod}), M)$ factor through $\varphi$.
\smallskip

\noindent {\bf (2)} If $\pinf(\Delta\mbox{\rm-mod})$ is contravariantly finite in $\Delta\mbox{\rm-mod}$, then $\Delta\mbox{\rm-mod}$ contains a unique \emph{basic} strong tilting module $T$ up to isomorphism $($that $T$ be \emph{basic} means that $T$ does not have any repeated indecomposable direct summands$)$.  Moreover, $\add(T)$ coincides with $\add(A(E))$, where $E$ is the minimal injective cogenerator in $\Delta\mbox{\rm-mod}$, and $A(E)$ is the unique right $\pinf(\Delta\mbox{\rm-mod})$-approximation of $E$ of minimal vector space dimension.
\end{theorem}
\smallskip

\noindent In this connection, we point out that the minimal approximations of the indecomposable injective modules may be decomposable.  
\medskip

\subsec*{3.B. Back to truncations of $\la = KQ$}
A major asset of the categories $\la_L\mbox{\rm-mod}$ is the fact that they always have strong tilting modules.  Our goal in this subsection is to explore the basic strong tilt $\latilde_L$, compare it with $\la_L$, and then explore the connections ``$\la_L$ versus $\latilde_L$" for $L \rightarrow \infty$.

\begin{theorem}
{\rm\cite[Theorem 4.1]{DHZ}}
For all $L \in \NN$,  the category $\pinf(\la_L\mbox{\rm-mod})$ is contravariantly finite in $\la_L\mbox{\rm-mod}$ and thus contains a strong tilting module.
\end{theorem}

\noindent 
Through Theorem \ref{tilt.summands}, we keep $L$ fixed and let $T = \bigoplus_{1 \le i \le n} T_i$ be the unique basic strong tilting object in $\la_L$-mod, where the $T_i$ are the indecomposable direct summands of $T$, one for each vertex of $Q$.  In fact, the vertices of $Q$ are in natural one-to-one correspondence with the $T_i$ and hence also with the vertices $\etil_i$ of the quiver $\Qtilde_L$ of $\latilde_L = \End_{\la_L}(T)^{\op}$, i.e., there are natural bijections $e_i \leftrightarrow T_i \leftrightarrow \etil_i$.  The latter assigns to $T_i$ the identity map on $T_i$ viewed as an element of $\End_{\la_L}(T) = \latilde_L^{\op}$; the correspondence $e_i \leftrightarrow T_i$  is detailed in Theorem \ref{tilt.summands}.

To obtain $T$ from $Q$ and $L$, we partition the vertices $e_1, \dots, e_n$ of $Q$ into precyclic and non-precyclic vertices: A vertex $e_i$ is \emph{precyclic} if there is a path in $Q$ which starts in $e_i$ and ends on an oriented cycle of $Q$. The precyclic vertices play a distinguished role in the homology of $\la_L\text{-mod}$ for arbitrary $L$ in that they coincide with the vertices that give rise to simple left $\la_L$-modules of infinite projective dimension. We carry over the attributes ``precyclic", resp.~``non-precyclic" defined for $e_i$ to the simple module $S_i = S_i^{(L)} \in \la_L\text{-mod}$ and the direct summand $T_i = T_i^{(L)}$ that corresponds to the vertex $e_i$ under the mentioned correspondence.  

The following idempotent of $\la_L$ is crucial in determining $\pinf(\la_L\mbox{\rm-mod})$-ap\-prox\-i\-ma\-tions, and hence in constructing the $T_i$, from $Q$ and $L$: 
$$\eps \ = \ \sum_{e_i \ \text{non-precyclic}} e_i\, .$$
Note that, for any left $\la_L$-module $M$, the subspace $\eps M$ is actually a submodule.
\smallskip

\begin{theorem}  \label{tilt.summands}
{\rm\cite[Theorems 4.2, 5.3]{DHZ}}  Let $L \in \NN$.
Given any finitely generated left $\la_L$-module $M$, the minimal right $\pinf(\la_L\mbox{\rm-mod})$-approximation of $M$ is the map
$$P(M)/\eps \Ker(f) \longrightarrow M \ \ \text{induced by the}\ \la_L\text{-projective cover}\  \ f: P(M) \rightarrow M.$$

\noindent The basic strong tilting module $T$ {\rm(see Theorem \ref{strong.tilting}(2))} can be explicitly constructed from $Q$ and $L$. Structure of the $T_i$: For any precyclic $e_i$, we have $T_i = \la_L e_i/ \eps J(\la_L) e_i$; in particular, all simple composition factors of $T_i$ are precyclic. When $e_i$ is non-precyclic, $T_i$ is the (indecomposable) minimal approximation of the injective envelope $E_i$ of $S_i$; in that case, $\eps T_i \cong \eps E_i$, $\ttop T_i \cong \ttop E_i$, and $T_i/\eps T_i$ is a direct sum of copies of precyclic $T_j$, one for each copy of $S_j$ in the top of $T_i$.  

In particular, $T_1, \dots, T_n$ have tree graphs and are  determined up to isomorphism by these graphs.
\end{theorem} 

\noindent 
We illustrate the preceding theorem with an example.

\begin{example}  \label{tilt.illustration}
Let $Q$ be the quiver below, $\la = KQ$, and $L=3$.
$$\xymatrixrowsep{1.5pc} \xymatrixcolsep{3.5pc}
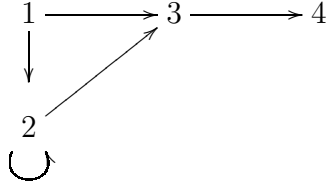
\xymatrix{
1 \ar[d] \ar[r] &3 \ar[r] &4  \\
2 \ar@(dl,dr) \ar[ur]  \vphantom{\biggl(}
}$$

\noindent In this case, the vertices $e_1, e_2$ are precyclic, $e_3, e_4$ non-precyclic, and the indecomposable direct summands of the basic strong tilting module $T \in \la_3\text{-mod}$ have graphs as follows:
$$\xymatrixrowsep{1.5pc} \xymatrixcolsep{0.4pc}
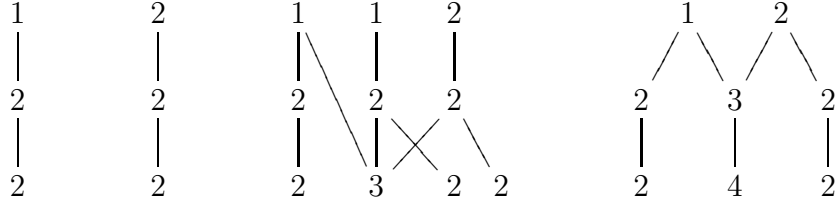
\xymatrix{
1 \edge[d] &&&&2 \edge[d] &&&&1 \edge[d] \edge[ddrr] &&1 \edge[d] &&2 \edge[d] & &&&& &1 \edge[dl] \edge[dr] &&2 \edge[dl] \edge[dr]  \\
2 \edge[d] &&&&2 \edge[d]  &&&&2 \edge[d] &&2 \edge[d] \edge[drr] &&2 \edge[dll] \edge[dr] & &&&&2 \edge[d] &&3 \edge[d] &&2 \edge[d]  \\
2 &&&&2 &&&&2 &&3 &&2 &2 &&&&2 &&4 &&2
}$$
Here the first two graphs, for $T_1 = \la e_1/\eps Je_1$ and $T_2 = \la e_2/\eps Je_2$, represent the minimal $\pinf(\lalmod)$-approximations of the simples $S_1$ and $S_2$.  The remaining two graphs depict  the minimal $\pinf(\lalmod)$-approximations of $E(S_3)$ and $E(S_4)$, as is verified via Theorem \ref{tilt.summands}.
\end{example}

\noindent 
The first of the main results of this section (Theorem \ref{B.A&A1}) targets the intrinsic structure of $\latilde_L$. Since a full description of a $K$-basis would involve a significant amount of notation, we only state availability of such a description.  As  can be gleaned from Example \ref{BAA.example} below, the basic $K$-algebra $\latilde_L \cong  K \Qtilde_L/ \widetilde{I}_L$, where $\Qtilde_L$ is the quiver of $\latilde_L$ and $\widetilde{I}_L \subseteq K \Qtilde_L$ a suitable admissible ideal, is no longer a monomial algebra in general; in Example \ref{BAA.example}, for instance, there is no choice of arrows in $\End_{\la_L}(T)$ that resolves all non-monomial relations displayed in the graph of the indecomposable projective left $\End_{\la_L}(T)$-module $\End_{\la_L}(T)\, \etil_2 = \etil_2\, \latilde_L$.

\begin{theorem}  \label{B.A&A1}
\cite{B+A&A} {\rm(loosely phrased)}
Let $Q$ be any quiver and $L \ge 2$.  Bases for the $\Hom$-spaces $\Hom_{\la_L}(T_i, T_j)$ can be explicitly obtained from $Q$ and $L$.  These bases are made up of five essentially different types of maps, next to certain compositions of these prototypical maps.

Consequence: A presentation of $\latilde_L$ via quiver and relations is readily available from $Q$ and $L$.
\end{theorem}

\noindent 
Combined with a hefty dose of combinatorics, Theorem \ref{B.A&A1} allows for comparisons of $\la_L$ with its strong tilt. In Theorems \ref{B.A&A2} and \ref{ratios}, we compare the Loewy lengths of the pair $\bigl(\la_L,\latilde_L\bigr)$ and explore the development of the ratios of these invariants as $L \rightarrow \infty$.

\begin{theorem}  \label{B.A&A2}
\cite{B+A&A}
As $Q$ ranges over arbitrary quivers, the Loewy length of $\latilde_L$ is uniformly bounded in terms of that of $\la_L$. Namely,
$$\bigl(\text{\rm Loewy length of} \ \la_L\bigr)  \ \le \ \bigl(\text{\rm Loewy length of} \ \latilde_L\bigr) \ \le \  3 \cdot  \bigl(\text{\rm Loewy length of}\ \la_L\bigr) - 2$$
for all quivers $Q$ and all positive integers $L$. These bounds are sharp in the following sense: For each $L$ there exists a quiver $Q$ such that the Loewy length of $\latilde_L$ attains the upper bound; the same is true for the lower bound. 
\end{theorem}

\begin{theorem}  \label{ratios}
{\rm\cite{asymptotic}}
Let $Q$ be any quiver.

\noindent{\bf(1)} The set 
$$A(Q) : =\biggl\{\,\text{accumulation points of}\ \biggl(\frac{\text{\rm Loewy length of} \ \latilde_L}{\text{\rm Loewy length of}\ \la_L}\biggr)_{L=2}^\infty \biggr\}$$
is  finite and can be computed from the quiver $Q$. 
\smallskip

\noindent{\bf (2)} The lim sup of the sequence of fractions
\begin{equation}  \label{ratiosequence}
\biggl(\frac{\text{\rm Loewy length of} \ \latilde_L}{\text{\rm Loewy length of}\ \la_L}\biggr)_{L=2}^\infty
\end{equation}
\smallskip
belongs to the interval $[1,2]$.
\smallskip

\end{theorem}

\noindent 
Of course, due to finiteness of $A(Q)$, the lim sup of the sequence \eqref{ratiosequence} is the maximum of the set $A(Q)$.  The reason for finiteness of $A(Q)$ lies in the fact that the Loewy length of $\latilde_L$ is completely determined by the (tree) graphs of the two uppermost radical layers of the indecomposable injective left $\End_\la(T)$-modules $E(S_i)$; for $L \gg 0$, these layers repeat periodically.
\medskip

\noindent The following example illustrates the construction of a quiver presentation of $\latilde_L$ from $Q$ and $L$ in the special case where $L=3$. The displayed class of quivers $Q^{(L)}$ also attests to sharpness of the upper bound for the Loewy length of $\latilde_L$ given in Theorem \ref{B.A&A2}.  We note that the lower bound is obvious, as is its sharpness, and hence we do not address it with (ubiquitous) examples.
  
\begin{examples}  \label{BAA.example}
For $L \ge 3$, let $Q^{(L)}$ be the quiver below, and let $\la_L := (KQ^{(L)})_L$ be the truncation of Loewy length $L$ of the path algebra $KQ^{(L)}$.  (For $L=3$, we re-encounter the quiver of Example \ref{tilt.illustration}.)
$$\xymatrixrowsep{1.5pc} \xymatrixcolsep{3.5pc}
\xymatrix{
1 \ar[d] \ar[r] &3 \ar[r] &4 \ar[r] &\ar@{}[r]|{\cdots} &\ar[r] &L+1  \\
2 \ar@(dl,dr) \ar[ur] \vphantom{\biggl(}
}$$
Then the Loewy length of $\latilde_L$ is $3L-2$ \cite{B+A&A}.

We back this up for $L=3$, based on the structure of the indecomposable direct summands of the strong tilting module $T \in \latilde_3\text{-mod}$ displayed in Example \ref{tilt.illustration}.  Theorem \ref{B.A&A1} allows us to compute quiver and relations of $\End_\la(T) = \latilde_3^{\op}$.  The quiver is as follows: 
$$\xymatrixrowsep{2pc} \xymatrixcolsep{3.5pc}
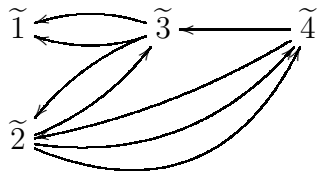
\xymatrix{
\onetil &\threetil \ar@/_/[l] \ar@/^/[l] \ar@/_/[dl] &\fourtil \ar[l] \ar@/^/[dll]  \\
\twotil \ar@/_/[ur] \ar@/_1.5pc/[urr] \ar@/_2.5pc/[urr]
}$$
\medskip

\noindent Clearly, the vertex $\etil_1$ is a sink, whence $\etil_1 \latilde_3$ is simple. For best visualization of the projective right $\latilde_3$-module $\tilde{e}_2\, \latilde_3$, alias the indecomposable projective left $\End_\la(T)$-module $\End_\la(T)\, \etil_2$, we present a graph: 
$$\xymatrixrowsep{1.5pc} \xymatrixcolsep{1.25pc}
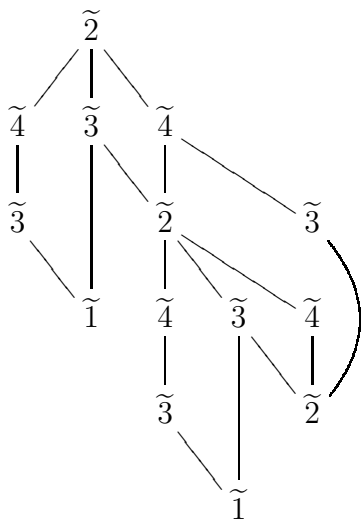
\xymatrix{
&\twotil \edge[dl] \edge[d] \edge[dr]  \\
\fourtil \edge[d] &\threetil \edge[dd] \edge[dr] &\fourtil \edge[d] \edge[drr] \\
\threetil \edge[dr] &&\twotil \edge[d] \edge[dr] \edge[drr] &&\threetil \edge@/^1.5pc/[dd]  \\
&\onetil &\fourtil \edge[d] &\threetil \edge[dd] \edge[dr] &\fourtil \edge[d]  \\
&&\threetil \edge[dr] &&\twotil  \\
&&&\onetil
}$$
\noindent We read off that the Loewy length of $\tilde{e}_2 \latilde_3$ is $6$. Since, via the arrow from $\etil_3$ to $\etil_2$, the module $\tilde{e}_2 \latilde_3$ embeds into the radical of  $\tilde{e}_3 \latilde_3$, this brings the Loewy length of $\tilde{e}_3 \latilde_3$ up to $7 = 3\cdot 3 - 2$.  Indeed, one obtains a nonzero composition of maps 
$T_3 \rightarrow T_2 \rightarrow T_4 \rightarrow T_2 \rightarrow T_4 \rightarrow T_3 \rightarrow T_1$ in $\End_\la(T)$, for instance. 

Finally, one checks that, for any $L \ge 2$, the sequence $\biggl(\dfrac{\text{\rm Loewy length of} \ (\widetilde{KQ^{(L)}})_{L'}}{\text{\rm Loewy length of}\ (KQ^{(L)})_{L'}}\biggr)_{L'=2}^\infty$ converges to $2$.
\end{examples}
\smallskip

\noindent The sequence  \eqref{ratiosequence} holds a great deal of information about the quiver $Q$, in particular about the lengths and connectivity of its oriented cycles.  It can be computed from $Q$. In the next example, we do so. 

\begin{example}  \label{example.tilt.sequence}
Take $Q$ to be
$$\xymatrixrowsep{1pc} \xymatrixcolsep{0.8pc}
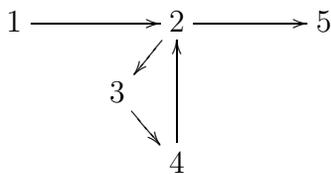
\xymatrix{
1 \ar[rrr] &&&2 \ar[dl] \ar[rrr] &&&5  \\
&&3 \ar[dr]  \\
&&&4 \ar[uu]
}$$
Here the Loewy length of $\la_L$ is $L$, since $Q$ has an oriented cycle. Denoting by $\widetilde{L}$ the Loewy length of $\latilde_L$, one finds:  If $L$ is at least $2$ and not divisible by $3$, then $\widetilde{L} = L + 1$; if, on the other hand, $L$ is divisible by $3$, then $\widetilde{L} = 4L/3$. Thus, the sequence $\bigl(\widetilde{L}/L\bigr)_{L = 2}^\infty$ is 
$$\biggl(\frac{3}{2}, \, \frac{4}{3}, \,  \frac{5}{4}, \, \frac{6}{5}, \, \frac{4}{3},\,  \frac{8}{7}, \, \frac{9}{8}, \, \frac{4}{3}, \dots\biggr),$$
and $A(Q) = \{1, 4/3\}$.

\end{example}
\medskip

\noindent We single out two of the most obvious questions left open so far. 

\begin{problems} 
 
$\bullet$ Is $A(Q)$ always a subset of $\QQ$?
 \smallskip

\noindent $\bullet$ Do arbitrary rational numbers in $[1,2]$ arise as maxima of sets $A(Q)$ for suitable $Q$?
\end{problems}

\noindent Further applications of Theorem \ref{B.A&A1}  follow in the next section.

\sectionnew{\bf Representation type}

\noindent Again let $\la = KQ$ and $\la_L$ the truncation of $\la$ that annihilates the paths of length $\ge L$.  There are two major unresolved problems concerning the representation type of these algebras:
\begin{itemize}
\item Describe the pairs $(Q, L)$ with the property that $\la_L$ \ has finite / tame infinite / wild representation type.
\item Describe the pairs $(Q, L)$ such that the strong tilt  $\latilde_L$ has finite / tame infinite / wild type.
\end{itemize}
The exploration of these problems is just beginning. 
\medskip

An obvious sufficient condition for any path algebra modulo relations to have infinite representation type is the occurrence of double arrows in its quiver;  triple arrows obviously lead to wild representation type.  These are criteria that are easy to check.

Obviously, ``most" quivers $Q$ with oriented cycles lead to algebras $\la_L$ and $\latilde_L$ of wild representation type for $L \gg 0$, prompting the question ``for which values of $L$ do switches occur?" 

Based on  $\Qtilde_2$ and the relations of $\latilde_2$, it is not difficult to construct quivers $Q$ with the property that $\la_2$ has finite type, whereas $\latilde_2$ has  any prescribed representation type. The accessibility of the case $L = 2$ is due to the fact that the separated quiver of $Q$ (e.g., \cite[p.350]{ARS}) allows to detect the representation type of any $(J^2 = 0)$-algebra by means of Gabriel's Theorem \cite{GabThm} and the supplements provided by Donovan-Freislich \cite{DoFr} and Nazarova \cite{Naz}.  

\begin{examples}  \label{finite.tilt.anything}
In each of the examples below, let $\la = KQ$ where $Q$ is the given quiver.  We use \cite[Theorem X.2.6]{ARS} to verify that $\la_2$ has finite type in all three cases.
\smallskip

{\bf(1)} 
Let $Q$ be the quiver $\xymatrixcolsep{2.5pc}
\xymatrix{
1 \ar@(dl,ul)^{\phantom{.}} \ar[r] &2
}$.  The separated quiver of $Q$ has the form
$$\xymatrixrowsep{1.5pc} \xymatrixcolsep{2pc}
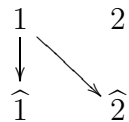
\xymatrix{
1 \ar[d] \ar[dr] &2  \\
\widehat{1} &\widehat{2}
}$$
which is a disjoint union of the Dynkin diagrams $A_2$ and $A_1$.  Consequently, $\la_2$ has finite representation type.

It is easily checked that $\latilde_2$ is a monomial algebra with quiver 
$$\xymatrixcolsep{3.5pc}
\xymatrix{
\onetil \ar@/^/[r] &\twotil \ar@/^/[l]
}$$
Since this quiver is an oriented cycle, $\latilde_2$ is a Nakayama algebra and thus has finite representation type (e.g., \cite[Theorems V.3.2, V.3.5]{AsSiSk}).
\smallskip

{\bf(2)}  Let $Q$ be the quiver $\xymatrixcolsep{2.5pc}
\xymatrix{
1 \ar@(dl,ul)^{\phantom{.}} \ar[r] &3 &2 \ar[l] \ar@(ur,dr)_{\phantom{m}}
}\;$.  This time the separated quiver of $Q$ has the form
$$\xymatrixrowsep{1.5pc} \xymatrixcolsep{2pc}
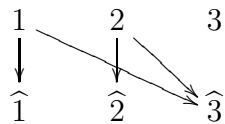
\xymatrix{
1 \ar[d] \ar[drr] &2 \ar[d] \ar[dr] &3  \\
\widehat{1} &\widehat{2} &\widehat{3}
}$$
which is again a disjoint union of simply laced Dynkin diagrams.  Thus again $\la_2$ has finite representation type.
\smallskip

The algebra $\latilde_2^{\op} = \End_{\la_2}(T)$ turns out to be a monomial algebra with quiver
$$\xymatrixcolsep{3.5pc}
\xymatrix{
\onetil \ar@/^/[r] &\threetil \ar@/^/[l]  \ar@/_/[r] &\twotil \ar@/_/[l] 
}$$
and indecomposable left projectives with graphs
$$\xymatrixrowsep{1pc} \xymatrixcolsep{0.4pc}
\xymatrix{
\onetil \edge[d] &&&\twotil \edge[d] &&& &\threetil \edge[dl] \edge[dr]  \\
\threetil \edge[d] &&&\threetil \edge[d] &&&\onetil \edge[d] &&\twotil \edge[d]  \\
\onetil &&&\twotil &&&\threetil \edge[d] &&\threetil \edge[d]  \\
&&& &&&\onetil &&\twotil
}$$
In particular, $\latilde_2^{\op}$ is a special biserial algebra, so it has tame representation type \cite[Corollary 2.4]{WW}.  Clearly, this algebra has indecomposable string modules of arbitrarily high dimension with graphs of the form
$$\xymatrixrowsep{1pc} \xymatrixcolsep{0.25pc}
\xymatrix{
&&\threetil \edge[dl] \edge[dr] &&&&\threetil \edge[dl] \edge[dr] && &&& &&\threetil \edge[dl] \edge[dr]  \\
&\onetil \edge[dl] &&\twotil \edge[dr] &&\onetil \edge[dl] &&\twotil \edge[dr] &&&&&\onetil \edge[dl] &&\twotil \edge[dr]  \\
\threetil &&&&\threetil &&&&\threetil &\cdots &\cdots &\threetil &&&&\threetil
}$$
Therefore $\latilde_2^{\op}$, and so also $\latilde_2$, has tame infinite representation type.
\smallskip

{\bf(3)} Finally, let $Q$ be the quiver $\xymatrixrowsep{0.5pc}\xymatrixcolsep{2.5pc}
\vcenter{\xymatrix{
\text{\phantom{.}}  \\
2 &1 \ar@(ul,ur) \ar[r] \ar[l] &3
}}\;$.  The corresponding separated quiver has the form
$$\xymatrixrowsep{1.5pc} \xymatrixcolsep{2pc}
\xymatrix{
1 \ar[d] \ar[dr] \ar[drr] &2 &3  \\
\widehat{1} &\widehat{2} &\widehat{3}
}$$
and we once again conclude that $\la_2$ has finite representation type.
\smallskip

This time, $\latilde_2^{\op}$ is a monomial algebra with quiver
$$\xymatrixcolsep{3.5pc}
\xymatrix{
\twotil \ar@/^/[r] &\onetil \ar@/^/[l]  \ar@/_/[r] &\threetil \ar@/_/[l] 
}$$
and indecomposable left projectives
$$\xymatrixrowsep{1pc} \xymatrixcolsep{0.4pc}
\xymatrix{
&&\onetil \edge[dl] \edge[dr] && &&& &\twotil \edge[d] & &&& &\threetil \edge[d]  \\
&\twotil \edge[d] &&\threetil \edge[d] & &&& &\onetil \edge[dl] \edge[dr] & &&& &\onetil \edge[dl] \edge[dr]  \\
&\onetil \edge[dl] \edge[d] &&\onetil \edge[d] \edge[dr] & &&&\twotil &&\threetil &&&\twotil &&\threetil  \\
\twotil &\threetil &&\twotil &\threetil
}$$
There is a $2$-parameter family of pairwise non-isomorphic indecomposable modules which share the graph 
$$\xymatrixrowsep{1pc} \xymatrixcolsep{0.4pc}
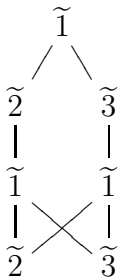
\xymatrix{
&\onetil \edge[dl] \edge[dr]  \\
\twotil \edge[d] &&\threetil \edge[d]  \\
\onetil \edge[d] \edge[drr] &&\onetil \edge[d] \edge[dll]  \\
\twotil &&\threetil
}$$
\smallskip

\noindent Indeed, if the arrows in the quiver of $\latilde_2^{\op}$ are named as follows
$$\xymatrixcolsep{3.5pc}
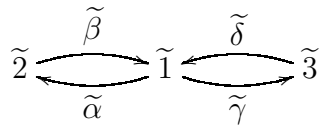
\xymatrix{
\twotil \ar@/^/[r]^{\widetilde{\beta}} &\onetil \ar@/^/[l]^{\widetilde{\alpha}}  \ar@/_/[r]_{\widetilde{\gamma}} &\threetil \ar@/_/[l]_{\widetilde{\delta}} 
}$$
we obtain a family of $\latilde_2^{\op}$-modules $M_{k,l} = \latilde_2^{\op}\etil_1/ U_{k,l}$, where $U_{k.l}$ is generated by the differences $\widetilde{\alpha} \widetilde{\beta} \widetilde{\alpha}  - k \cdot \widetilde{\alpha} \widetilde{\delta} \widetilde {\gamma}$ and $\widetilde{\gamma} \widetilde{\delta} \widetilde{\gamma} - l \cdot \widetilde{\gamma} \widetilde{\beta} \widetilde{\alpha}$.  To check distinctness of the isomorphism types, suppose there is an isomorphism $\phi: M_{k,l} \rightarrow M_{k',l'}$.  If $x := \etil_1 + U_{k,l}$ and $x' := \etil_1 + U_{k',l'}$, then $\phi(x)$ must be a top element of $M_{k',l'}$ of the form $ax' + b \widetilde{\beta} \widetilde{\alpha} x' + c \widetilde{\delta} \widetilde{\gamma} x'$  with $a \in K^*$, $b, c \in K$.  Since $k \widetilde{\alpha} \widetilde{\delta} \widetilde {\gamma} \phi(x) = \widetilde{\alpha} \widetilde{\beta} \widetilde{\alpha} \phi(x) = a \widetilde{\alpha} \widetilde{\beta} \widetilde{\alpha} x' = a k' \widetilde{\alpha} \widetilde{\delta} \widetilde {\gamma} x' = k' \widetilde{\alpha} \widetilde{\delta} \widetilde {\gamma} \phi(x)$, we obtain $k = k'$, and similarly $l=l'$.

Therefore $\latilde_2^{\op}$, and thus $\latilde_2$, has wild representation type in this case.
\end{examples}

\noindent 
For $L = 2$, any quiver $Q$ with the property that $\la_2$ has finite representation type gives rise to a strong tilt $\latilde_2$ whose quiver is free of double arrows.  

However, in general the multiplicity of equidirected arrows between two given vertices of the quiver of $\latilde_L$ may grow exponentially as a function of $L$, even when the original quiver $Q$ does not have double arrows. 
We conclude with a simple example. 

\begin{example}
Let $Q$ be the quiver below, and $\la = KQ$.
$$\xymatrixrowsep{0.8pc} \xymatrixcolsep{3.5pc}
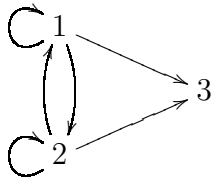
\xymatrix{
1 \ar@(dl,ul) \ar@/^/[dd]  \ar[dr]   \\
&3  \\
2 \ar@(dl,ul)  \ar@/^/[uu]  \ar[ur] 
}$$
Then the number of arrows from $\etil_3$ to $\etil_1$ in the quiver of $\latilde_L$ is $ \ge 2^{L-2}$ for $L \ge 2$.
\end{example}

\vfill\eject


\end{document}